\documentclass[12pt,reqno]{amsart}
\usepackage{mathtools,footnote}
\usepackage{amsmath,amsthm,amsfonts,amssymb}
\usepackage{verbatim,graphicx}

\renewcommand{\epsilon}{\varepsilon}
\usepackage[usenames, dvipsnames]{color}
\usepackage{a4wide}
\usepackage{url}

\newcommand{\im}{\textnormal{Im}}
\newcommand{\re}{\textnormal{Re}}

{\theoremstyle{definition}
}

\title[Prime numbers and 
Riemann zeros]{A computational history of prime numbers and Riemann zeros}
\author[P. Moree]{Pieter Moree}
\author[I. Petrykiewicz]{Izabela Petrykiewicz}
\author[A. Sedunova]{Alisa Sedunova}
\date{\today}

\address[P. Moree, A. Sedunova]{%
Max-Planck-Institut f\"ur Mathematik\\
Vivatsgasse 7\\D-53111 Bonn\\Germany.}
\email{moree@mpim-bonn.mpg.de}
\email{alisa.sedunova@phystech.edu}

\address[I. Petrykiewicz]{}
\email{ipetrykiewicz@gmail.com}

\date{}


\begin{document}

\maketitle 
\begin{abstract}
We give an informal survey of the historical development of computations related to prime number distribution 
and zeros of the Riemann zeta 
function\footnotemark.
\end{abstract}

\footnotetext{Caution! The authors are non-experts.}

The fundamental
quantity in the study of prime numbers\footnote{We follow the 
tradition to denote a prime number by the letter $p$.} 
is the \textbf{prime counting function} $\pi(x)$, which counts the number of primes not 
exceeding $x;$ 
in mathematical notation we have
$$\pi(x)=\sum_{p\le x}1.$$
\indent The first mathematicians to investigate the
growth of $\pi(x)$ had of course to start with collecting
data. They did this by painfully setting up tables of consecutive 
prime numbers, e.g., Kr\"uger in 1746 and 
Vega in 1797 (primes up to $100\,000$ and $400\,031$ respectively).
The most celebrated
of these prime table computers was Gauss. 
In 1791, when he was 14 years old, he noticed that as one gets to larger and larger numbers
the primes thin out, but that locally their distribution appears to be quite 
erratic. He based himself on a prime number table contained in 
a booklet with tables of logarithms he had received as a prize, 
and went on to conjecture that the ``probability that 
an arbitrary integer $n$ is actually a prime number should equal $1/\log n$''.
Thus Gauss conjectured that
$$\pi(x)\approx \sum_{2\le n\le x}\frac{1}{\log n}\approx {\rm Li}(x),$$ with
$${\rm Li}(x)=\int_2^x \frac{dt}{\log t},$$
the \textbf{logarithmic integral}\footnote{An ever recurring theme in analytic number theory
is approximating a sum by an integral (Section \ref{mac}).}.
Since by partial integration it
is easily seen that Li$(x)\sim x/\log x$, the conjecture of Gauss implies that asymptotically
$$\pi(x) \sim \frac{x}{\log x},$$
a conjecture that was proved much later, in 1896, by 
Hadamard\footnote{Jacques Salomon Hadamard (1865 Versailles, France -- 1963 in Paris, France), French mathematician, 
professor at the University of Bordeaux, Coll\`{e}ge de France, \'{E}cole Polytechnique and \'{E}cole Centrale.} 
and 
de la Vall\'{e}e-Poussin\footnote{Charles Jean Gustave Nicolas Baron de la Vall\'{e}e-Poussin (1866 Leuven, Belgium -- 1962 Leuven, Belgium), Belgian mathematician,
professor at Catholic University of Leuven, Coll\`{e}ge de France and Sorbonne.} independently.
This asymptotic for $\pi(x)$ is called the
\textbf{Prime Number Theorem} (PNT).\\
\indent Gauss kept a life long interest in primes and what
he did was to count primes in blocks of $1\,000$ (a Chiliade). 
As he wrote 
in a letter to Bessel, he would
use an idle quarter of hour here and there to deal with a further block. By 
the end of his life he would extend the tables up to $3\,000\,000$.
After Gauss, number theorists kept extending the existing prime number
tables. Thus in 1856 
Crelle\footnote{August Leopold Crelle (1780 Eichwerder, Germany -- 1855 Berlin, Germany), German mathematician, founder of \textit{Journal f\"{u}r die reine und angewandte Mathematik}.}
published a 
table of primes up to $6\,000\,000$, and a few years later 
Dase\footnote{Johann Martin Zacharias Dase (1824 Hamburg, Germany -- 1861 Hamburg, Germany), German arithmetician, having great calculating skills, but little mathematical knowledge.} 
extended this to $9\,000\,000$. The most impressive feat in this regard is
due to 
Kulik\footnote{Jakob Philipp Kulik (1793 Lemberg, Austrian Empire -- 1863 Prague, Bohemia), Polish-Austrian mathematician, professor at the Charles University of Prague.},
who spent 20 years preparing a 
factor table of the numbers coprime to $30$
 up to $1\,000\,330\,200$ (he did so in eight 
manuscript volumes, totalling $4\,212$ pages).\\
\indent The holy grail in computational prime number theory is to find sharp 
estimates of $\pi(x)$. These
estimates should be in terms of elementary functions. 

An early attempt is by 
Legendre\footnote{Adrien-Marie Legendre (1752 Paris, France -- 1833 Paris, France), French mathematician and author of an influential number theory book.}, who claimed (1808) that
$x/(\log x-1.0836)$ should approximate $\pi(x)$ well. We now
know that this is a reasonable estimate (the estimate
$x/(\log x-1)$ is actually better).
A much more recent and rigorous example is provided by
the estimates
$$\frac{x}{\log x}\Big(1+\frac{1}{2\log x}\Big)<\pi(x)<
\frac{x}{\log x}\Big(1+\frac{3}{2\log x}\Big), \quad x\ge 59.$$
due to Rosser and Schoenfeld 
\cite{MR0137689}. Some further 
examples can be found in Section \ref{expli}. The reason why
sharp estimates of $\pi(x)$ 
and of related prime counting functions
are so important is that many problems
in number theory use them as input. There are plenty
of number theoretical problems where one comes to a solution only on assuming that
a sharp estimate for $\pi(x)$ is available, an estimate we cannot currently
prove, but which we could if we knew that the {\bf Riemann Hypothesis} (RH) holds
true (we will come back to this shortly). Under RH it can be shown
that for every $x>2\,657$ we have
\begin{equation}
\label{schoen1969}
|\pi(x)-{\rm Li}(x)|<\frac{1}{8\pi} \sqrt{x}\log x.
\end{equation}
This is a sharp inequality as the 
estimate $\pi(x)={\rm Li}(x)+O(\sqrt{x}/\log x)$ does not
hold \cite{Schm}. The latter result says that the 
primes do behave irregularly to some extent.
The importance of the RH is that if true it would imply that we
can very well approximate $\pi(x)$ by the simple function Li$(x)$, which
makes proving results involving primes in general much easier.\\
\indent The sharpest estimates to date for $\pi(x)$ are obtained 
by using properties
of the so-called {\bf Riemann zeta function}, defined by
\begin{equation}
\label{zetadefi}
 \zeta(s) = \sum_{n=1}^\infty \frac{1}{n^s},
\end{equation}
with $s$ a complex number having real part $\re s>1$, see \cite{R}\footnote{The imaginary part
of $s$ will be denoted by $\im s$.}.
The function converges for all complex numbers $s$ such that 
$\re s > 1$. 
In 1859 a renowned G\"{o}ttingen professor, 
Riemann\footnote{Georg Friedrich Bernhard Riemann (1826 Breselenz, Hanover -- 1866 Selasca, Italy), German mathematician, student of Gauss, professor at G\"{o}ttingen University.},
published
``\"{U}ber die {A}nzahl der {P}rimzahlen unter einer gegebenen {G}r\"{o}sse''.
This paper\footnote{This only 9-pages-long paper, it is the only published work of Riemann on 
number theory.} is without doubt the most important paper ever written in analytic number theory; indeed it is foundational as Riemann makes essential use
of $s$ being a complex variable, whereas a century 
earlier Euler\footnote{Leonhard Euler (1707 Basel, Switzerland -- 1783 St. Petersburg, Russian Empire), extremely prolific Swiss mathematician, student of Johann Bernoulli.} only considered $\zeta(s)$
for real values of $s$\footnote{For more on Euler's work on $\zeta$, see
Ayoub \cite{MR0360116}.}. That Riemann considered his zeta function as an analytic
function comes as no surprise as the development of complex analysis was one of his 
central preoccupations.
He brought the study of $\pi(x)$ to a completely new level, but 
actually proved little as the tool he used, complex function theory, did not have 
on a firm theoretical foundation at the time. It took 
about 40 years and a lot of preliminary work, mainly
by Cahen, Halphen and Phragm\'en (see, e.g., \cite{N} for 
more details), before the PNT 
could be finally proved using methods of complex function 
theory. A tremendous amount of work was carried out by
Landau\footnote{Edmund Georg Hermann Landau (1877 Berlin, Germany -- 1938 Berlin, Germany), German mathematician, student of Frobenius, professor at G\"{o}ttingen University.}, 
who went meticulously through all earlier relevant work on this subject, checked its correctness, simplified it\footnote{He gave a much simpler proof of the PNT, for
example.}
and wrote a standard work on prime number theory in 1909 \cite{La, MR0068565}\footnote{H. Montgomary and R. Vaughan, both non-Germans, studying the original German version, were surprised to learn about a very strong mathematician called Verfasser they had never heard of (Verfasser means
author...). }. 
He himself proved many important results as well, such 
as the prime ideal theorem (see Section \ref{zz}).\\
\indent A lot of effort was put into proving results about prime numbers 
without using complex analysis.  The most celebrated results were obtained
by Selberg\footnote{Atle Selberg (1917
Langesund, Norway -- 2007
Princeton, the US), 
Norwegian mathematician, professor at Princeton University, recipient of 
the Fields Medal (1950) and the Abel Prize (2002).}
(and, more or less, independently) by Erd\H{o}s. They based
themselves on the identity, now called \textbf{Selberg's symmetry formula},
$$\sum_{p\le x}\log^2 p+\sum_{p,q\le x}\log p \cdot \log q=2x\log x+O(x),$$
and used it to obtain an \textbf{elementary proof} of the PNT. 
Until 1950 it was widely believed (e.g., by Landau) that no such elementary proof could be developed and so this result greatly impressed the contemporaries.
Later Selberg extended this combinatorial technique to show the PNT for 
primes in arithmetic progression, see \cite{MR0033306}. For a nice
survey see Diamond \cite{MR670132}. Unfortunately, the high hopes placed in
new insights coming from finding an elementary proof of PNT were
thwarted.\\
\indent The level
of insight that Riemann had reached was finally surpassed by Hardy and Littlewood in the 1920's. As far as his zeta function is concerned, Riemann was certainly more than half a century ahead of his time!\\
\indent Using partial integration it is easy to
deduce that for $\re s>0$ an analytic continuation
of $\zeta(s)$ is given by
\begin{equation}
\label{simplecont}
\zeta(s)=\frac{s}{s-1}-s\int_1^{\infty}\frac{t-\lfloor{t\rfloor}}{t^{1+s}}dt,
\end{equation}
with $\lfloor t \rfloor$ being the floor function. Another analytic continuation is obtained
on noting that for $\re s>1$ we have
\begin{equation}
\label{simplecont2}
(1-2^{1-s})\zeta(s)=\sum_{n=1}^{\infty}(-1)^n n^{-s}.
\end{equation}
Since the right hand side actually converges for $\re s>0$, the identity furnishes an 
analytic continuation for all $\re s>0$.\\
\indent In his paper, Riemann showed that the zeta function 
actually has an analytic continuation to the 
{\textbf {whole}} complex plane, except for a simple pole at $s=1$,
and that it vanishes at all negative even integers. The \textbf{trivial zeros} of $\zeta(s)$ are the ones at negative even integers, the \textbf{non-trivial zeros} come from complex numbers $s=\sigma+it$ with $0 \leq \sigma \leq 1$.
The zeta function satisfies the {\bf functional equation}
\begin{equation}
\label{functionalvg}
\Gamma\left(s/2\right)\pi^{-s/2} \zeta(s)=\Gamma\left((1-s)/2\right)\pi^{(s-1)/2} \zeta(1-s),
\end{equation}
where $\Gamma(z)=\int_{0}^\infty t^{z-1}e^{-t}dt$ for $\re z>0$
denotes the \textbf{Gamma function} (a continuous extension of the \textbf{factorial 
function}\footnote{If $n$ is a non negative integer, then, e.g., $\Gamma(n+1)=n!$. For this reason some authors write $\Gamma(z+1)$ instead
of $\Gamma(z)$.}). This function is never equal to zero, and is holomorphic
everywhere except at the points $0,-1,-2,\ldots,$ where it has
simple poles.\\
\indent It turns out that $\Gamma(s/2)\pi^{-s/2}\zeta(s)$ has
a simple pole at both $s=0$ and $s=1$. This suggests 
that we should multiply
it by $s(s-1)$. In this way we obtain the
\textbf{Riemann} $\mathbf{\xi}$\textbf{-function}.
It is an \textbf{entire function} whose zeros $\rho$ are the non-trivial zeros of 
$\zeta( s )$. 
Note that we can rewrite (\ref{functionalvg}) as
$$\xi( s ) = \xi( 1-s ).$$
We have the following \textbf{Hadamard factorisation} of
$\xi$:
$$\xi(s)=s(s-1)\Gamma(s/2)\pi^{-s/2}\zeta(s)=\xi(0)\prod_{\rho}
\left(1-s/\rho\right)
e^{s/\rho},$$
where here 
(and in the sequel) the zeros 
$\rho$ are counted with 
their own multiplicities, e.g., a double zero is counted twice.\\
\indent Riemann gave two beautiful
proofs of the functional equation; one is related to the theory of modular forms and
makes use of the transformation property $\Omega(x)=x^{-1/2}\Omega(x^{-1})$, valid for
all positive $x$, with
$$\Omega(x)=\sum_{n=-\infty}^{\infty}e^{-n^2\pi x},$$
while the second proof uses the integral representation 
(\ref{eq:ZetaIntegral}) from Section \ref{RSformula}.
Note that if $\rho$ is a non-trivial zero, then so is $1-\rho$ 
by the functional equation. 
Moreover, since $\overline{\zeta(s)}=\zeta(\overline{s})$, we deduce that $\overline{\rho}$ and $1 - \overline{\rho}$ are also zeros.
Thus the zeros are symmetrically arranged about the {\bf half line} (also 
called the \textbf{critical line}) given by $\re s=1/2$ and also about the real axis. 
Therefore we often only calculate the zeros in the upper half plane.
Riemann shows that the non-trivial zeros lie in the {\bf critical
strip}\footnote{In this strip formula 
(\ref{zetadefi}) for $\zeta(s)$ does not apply, 
but formula (\ref{simplecont}) and (\ref{simplecont2}) do.}, the strip $0\le \re s\le 1$, and, furthermore,
are not real.
Moreover, he wrote that probably all its non-trivial zeros lie on the 
half line and continues: ``Certainly one would wish for a stricter proof here; I have meanwhile
temporarily put aside the search for this after some fleeting futile attempts,
as it appears unnecessary for the next objective of my investigation.''\footnote{\textit{Hiervon w\"{a}re allerdings ein strenger Beweis zu w\"{u}nschen; ich habe indess die Aufsuchung
desselben nach einigen fl\"{u}chtigen vergeblichen Versuchen vorl\"{a}ufig bei Seite gelassen, da er f\"{u}r den n\"{a}chsten Zweck meiner Untersuchung entbehrlich schien.}
\cite{R} Translation: David R. Wilkins.}
The \textbf{Riemann Hypothesis} states that all non-trivial zeros (the zeros in
the critical strip) of the zeta 
function are on the half line. We call the zeros of $\zeta(s)$ on the half line 
\textbf{critical zeros}. If in the rectangle 
defined by $0\le \re s\le 1$ and $0\le \im s\le T$ there are
only critical zeros, we say that \textbf{RH up to height $T$} 
holds true.\\
\indent Solving RH is one of the famous 
23 problems posed by Hilbert at the 1900 International Congress of 
Mathematicians in Paris. These problems held their fascination and influence 
on the developments through the twentieth century \cite{MR1930195}. Resolving RH is one of
seven Millennium Prize Problems that were stated by the Clay Mathematics Institute in 2000. A correct solution to any of the problems results in a $1\,000\,000$ dollar 
 prize being awarded by this institute \cite{MR2246251}.
Many mathematicians regard RH as the biggest
open problem in all of mathematics.\\
\indent The behaviour of
$\zeta(s)$ in the critical strip
is very closely related to the 
distributional properties of the primes. 
The uniqueness of prime factorization finds its analytic
counterpart in the identity
\begin{equation}
\label{prodie}
\zeta(s)=\prod_p (1-p^{-s})^{-1},~\re s>1,
\end{equation}
which was established by
Euler for real $s$. Indeed, to become an inhabitant of
the Zeta Zoo (Section \ref{zz}) an analytic function $f(s)$ needs
to have a factorization of the form
$\prod_p f_p(s).$ A formula of this type is
now called an \textbf{Euler product}.

Thus the fact
that $\zeta(s)$ tends to infinity if one approaches $s=1$ from the right over the real
axis\footnote{Indeed, the Riemann zeta
function has a simple pole with residue $1$ at $s=1$.} ensures 
by (\ref{prodie}) that there are infinitely many primes
$p$. This was discovered in 1737 by Euler who established the stronger 
result that
$\sum_{p}p^{-1}$ is unbounded\footnote{Gauss in 1796 conjectured, and Mertens proved
in 1874 the stronger result that $\sum_{p\le x}p^{-1}=\log \log x + C +o(1)$, with $C$ a constant.}.
It turns out that, in order to prove the PNT, it is enough to
show that there are no zeros on the line $\re s=1$.
The connection between the non-trivial zeros and the prime numbers is 
actually much closer than this.
Indeed, it can be shown that
\begin{equation}
\label{R1}
\pi(x)=R(x)-\sum_{\rho}R(x^{\rho}),
\end{equation}
where the sum is over all the non-trivial zeros $\rho$ (counted with multiplicities)
and 
\begin{equation}
\label{R2}
R(x)=\sum_{n=1}^{\infty}\frac{\mu(n)}{n}
{\rm Li}(x^{1/n})
\end{equation}
denotes the \textbf{Riemann function} and 
$\mu$ the M\"obius function (see Section \ref{merti}). Using
Gram's (1893) quickly converging power series 
$$R(x)=1+\sum_{n=1}^{\infty}\frac{1}{n\zeta(n+1)}\frac{(\log x)^n}{n!},$$
the Riemann 
function is computable. Two further expressions for $R(x)$ were found by
Ramanujan. He independently discovered the importance of $R(x)$ (around 1910) and developed
a prime number theory that, as Hardy phrased it, was "what the theory might be
if the zeta-function had no complex zeros".

\begin{table}
\footnotesize
 \begin{tabular}{ c | c |  c | c }
  $x$ & $\pi(x)$ & ${\rm Li}(x)-\pi(x)$  & $R(x)-\pi(x)$   \\
  \hline
$10^{8}$ & $5\,761\,455$ & $754$ & $97$ \\
$10^{9}$ & $50\,847\,534$ & $1\,701$ & $-79$ \\
$10^{10}$ & $455\,052\,511$ & $3\,104$ & $-1\,828$ \\
$10^{11}$ & $4\,118\,054\,813$ & $11\,588$ & $-2\,318$ \\
$10^{12}$ & $37\,607\,912\,018$ & $38\,263$ & $-1\,476$ \\
$10^{13}$ & $346\,065\,536\,839$ & $108\,971$ & $-5\,773$ \\
$10^{14}$ & $3\,204\,941\,750\,802$ & $314\,890$ & $-19\,200$ \\
$10^{15}$ & $29\,844\,570\,422\,667$ & $1\,052\,619$ & $73\,218$ \\
$10^{16}$ & $279\,238\,341\,033\,925$ & $3\,214\,632$ & $327\,052$ \\
$10^{17}$ & $2\,623\,557\,157\,654\,233$ & $7\,956\,589$ &  $-598\,255$ \\
$10^{18}$ & $24\,739\,954\,287\,740\,860$ & $21\,949\,555$ & $-3\,501\,366$ \\
\end{tabular}
\caption{The values of $\pi(x)$ compared with values of ${\rm Li}(x)$ and $R(x)$.}\label{table:pix}
\end{table}

\indent It turns out that, from a theoretical perspective, it is better
to work with certain weighted prime counting functions rather than with
$\pi(x)$ (that is, each prime $p$ is counted with a weight $w(p)$). The most well-known of these are $\vartheta(x)=\sum_{p\le x}\log p$
and $\psi(x)=\sum_{p^m\leq x}\log p$, where the sums are taken over all primes, respectively all prime powers 
less than $x$. The first function is known as the {\bf Chebyshev $\vartheta$-function}, 
the second one as the {\bf Chebyshev $\psi$-function}\footnote{It is
easy to see that $\psi(x)$ is the logarithm of the least common multiple of the integers 
$1,2,\dots,[x]$.}. 
One often sees the $\psi$-function defined as $\psi(x)=\sum_{n\le x}\Lambda(n)$, where
$\Lambda(n)$, the \textbf{von Mangoldt function}, equals $\log p$ if $n$ is a power of
a prime $p$ and zero otherwise. This function is more natural than it appears at first sight,
since logarithmic differentiation of the Euler identity 
(\ref{prodie}) yields
$-\zeta'/\zeta (s) =\sum_{n=1}^{\infty}\Lambda(n)n^{-s}$ for $\re s>1$. The quotient $\zeta'/\zeta$
has simple poles in the zeros of $\zeta$ and plays an important role in prime number theory.\\
\indent The complicated explicit formula (\ref{R1}) for $\pi(x)$
takes a much easier form for $\psi(x)$. This result
is called the ``explicit formula'' and
due to von Mangoldt\footnote{Hans Carl Friedrich von Mangoldt (1854 Weimar, Duchy of Saxe-Weimar-Eisenach -- 1925 Danzig-Langfuhr, Free City of Danzig), 
German mathematician, student of Kummer and Weierstrass, professor at Hanover University and Technical University of Aachen.}. 
He showed in 1895 that for $x>1$ and $x$ not
a prime power we have
\begin{equation}
\label{psipsi2}
 \psi(x) = x - \sum_{\zeta(\rho)=0} \frac{x^{\rho}}{\rho} - \frac{\zeta'(0)}{\zeta(0)}, 
\end{equation}
where the sum on the right-hand side is over all zeros (also the trivial ones!) of the Riemann zeta function. 
The explicit von Mangoldt formula was vastly
generalized by Andr\'e Weil. His explicit formula 
works for
a large class of test functions. The summation over
the roots involves the test function and this then
equals a sum over the primes involving the Fourier 
transform of the test function. On choosing a suitable test function (depending on the
problem one studies), Weil's explicit formula has 
actually become an important tool in computational prime number theory.\\
\indent Note that $|x^{\rho}|=x^{\re \rho}$ and consideration of the explicit formula suggests that
$\lim \sup \re \rho$ determines the error term.
Due to the existence of critical zeros, the minimum possible value
of the latter quantity is a half. We have $\lim \sup \re \rho =1/2$ if and
only if the RH holds true\footnote{Recall that if $\rho$ is a zero, so is
$1-\rho$.}. \\
\indent Note that if 
the $\sum 1/|\rho |$ were to be bounded, with 
$\rho$ ranging over all critical zeros, then it would
follow from (\ref{psipsi2}) that on RH 
we have $\psi(x)=x+O(\sqrt{x})$.
However, this sum does not converge, but a slightly weaker result is true, namely
$$\sum_{\im s \le T}\frac{x^{\rho}}{|\rho|}=O(x^{\alpha}\log^2 T),$$
with $\alpha$ any number such that there are no zeros $\rho$ with 
$\re \rho >\alpha$. Using the latter approximation with $\alpha=1/2$ 
it can be shown that the RH is equivalent with 
\begin{equation}
\label{psieee}
\psi(x)=x+O(x^{1/2}\log^2 x).
\end{equation}
By the same argument it can be unconditionally shown that
$\psi(x)=x+O({x}^{\kappa}\log^2 x)$
with $\kappa=\lim \sup \re \rho$.\\
\indent The analogue of (\ref{psieee}) for $\pi(x)$ is due to
von Koch, who showed in 1901 that the RH is equivalent to 
the formula
\begin{equation}
\label{pionRH}
\pi(x) = \textnormal{Li}(x)+O(E(x)), \quad E(x)=x^{1/2}\log x.
\end{equation}
As in the error term in (\ref{psieee}), the $1/2$ in the error term $E$ is 
directly related to the half line: if
we would have a zero off the half line, it always has a related zero
$s_0$ having $\re s_0>1/2$ and the estimate (\ref{pionRH}) with  
$E(x)=x^{\re s_0-\epsilon}$ is false. Indeed, more generally, the larger the 
area inside the critical strip where we can show there are no zeros, the smaller
we can take $E(x)$ in (\ref{pionRH}). Quite a bit of computational work was done
on determining an explicit zero free region that is as large
as possible, e.g., 
Kadiri \cite{MR2140161} showed that
$$\re s\ge 1-\frac{1}{R_0|\im s|},\quad |\im s|\ge 2,\quad R_0=5.69693,$$
is a zero free region. Already in 1899
de la Vall\'ee-Poussin had established the above result 
with different constants. He showed that
his zero free region leads to
$E(x)=\exp{(-c\sqrt{\log x})},$ with some positive constant $c$. It implies that for every
fixed $r>1$ we have $E(x)=x(\log x)^{-r}$. This function is far 
from behaving like $\sqrt{x}$. Indeed, it would be an astounding result if
somebody could prove that $E(x)=x^{\alpha}$, for some $\alpha<1$. 
Despite great effort by many number theorists, the above result 
of de la Vall\'ee-Poussin 
has not been much improved.\\
\indent As we have seen a prime number heuristic that works well is to
assume that $n$ is a prime with probability $1/\log n$. Riemann's research
leads one to replace this by a more accurate heuristic, namely that
the average value of $\Lambda(n)$ equals 1. This leads us to expect that, e.g., 
$\psi(x)\sim x$, which is the PNT in a different guise\footnote{Recall that
$\psi(x)=\sum_{n\le x}\Lambda(n)$.}. 
Indeed, by elementary arguments
it can be shown that the assertions
\begin{equation*}
\pi(x) \sim \frac{x}{\log x},\quad \vartheta(x)\sim x, \quad \psi(x)\sim x,
\end{equation*}
are all equivalent and so are three different guises of the PNT.\\
\indent More interesting is to consider
$\Pi(x)=\sum_{n\le x}\Lambda(n)/\log n$. Here our heuristic suggests that
\begin{equation}
\label{rielambda}
\Pi(x)=\sum_{k=1}^{\infty}\frac{\pi(x^{1/k})}{k}
\approx \sum_{n\le x}\frac{1}{\log n}\approx {\rm Li}(x).
\end{equation}
By M\"obius inversion we obtain from this the heuristic $\pi(x)\approx R(x)$, 
with $R(x)$ the Riemann function. This approximation is excellent, as shown in 
Table~\ref{table:pix}.\\
\indent From the papers Riemann
 left us, it transpires that his
computational skills 
were amazing and 
that he (being a perfectionist) kept a lot of his findings up his sleeve\footnote{Unfortunately a lot of his 
notes
were burnt after his death by his housekeeper.
A true treasure trove that went up in
cinders!}.
In a letter to Weierstrass from 1859 
Riemann mentioned that he discovered new expansions for $\zeta$. 
At the turn of the 20th century many people tried to find these expansions in the unpublished notes of Riemann in the G\"{o}ttingen library.
Meanwhile in England,
Hardy\footnote{Godfrey Harold Hardy (1877 Cranleigh, England~--~1947 Cambridge, England), English mathematician, professor at the University of Cambridge.}
in the early 1920's, together with his lifelong collaborator 
Littlewood,\footnote{John Edensor Littlewood (1885 Rochester, England -- 1977 Cambridge, England), British mathematician, professor at the University of Cambridge.} 
found ``an approximate functional equation'' \cite{HL1,HL2}. 
However, a librarian in 
G\"ottingen discovered that this approximate formula was already known to Riemann. 
Moreover, whereas Hardy and Littlewood had determined the dominating term only, Riemann had a method to estimate the remainder term.
In 1926, while in G\"{o}ttingen, 
Bessel-Hagen\footnote{Erich Paul Werner Bessel-Hagen (1898 Berlin, German Empire -- 1946 Bonn, Germany), German mathematician, student of Carath\'{e}odory, professor at the University of Bonn.} 
discovered a previously unknown approximate formula in Riemann's notes.
Siegel\footnote{Carl Ludwig Siegel (1896 Berlin, German Empire -- 1981 G\"{o}ttingen, West Germany), German mathematician, student of Landau, professor at the University of Frankfurt and G\"{o}ttingen.}
then analysed the notes with these two approximate formulas, 
working out the details,\footnote{\textit{Diese Gr\"{u}nde machten eine freie Bearbeitung des Riemannschen Fragmentes notwendig, wie sie im folgenden ausgef\"{u}hrt werden soll.} \cite{S}} 
and in 1932 published 
them. One of these
is now known as the \textbf{Riemann-Siegel formula}. 
This formula is still used in calculating the zeros of the Riemann zeta function, see Section \ref{RSformula} for more details.

Over time mathematicians calculated 
non-trivial zeros in 
the hope of disproving the Riemann Hypothesis (like Turing) or to find evidence for it.
At present, more than 100 billion zeros have been verified to be 
critical; 
for example, Gourdon in 2004 verified that the first $10^{13}$ zeros are
critical (see \cite{Go}); 
moreover, Odlyzko calculated 10 billion zeros near $10^{22}$ and showed that they all are critical, \cite{O}.
The first fifteen zeros are listed in Table~\ref{table:zeros}.
The number of non-trivial zeros computed by the various dramatis personae 
can be found in Table \ref{table:records}.

\begin{table}
 \begin{tabular}{ c | c ||  c | c ||  c | c  }
  $ n $ & zero &  $ n $ & zero &  $ n $ & zero \\
  \hline
  1 & 14.1347 & 6 & 37.5862 & 11 & 52.9703 \\
  2 & 21.0220 & 7 & 40.9187 & 12 & 56.4462 \\
  3 & 25.0109 & 8 & 43.3271 & 13 & 59.3470 \\
  4 & 30.4249 & 9 & 48.0052 & 14 & 60.8318 \\
  5 & 32.9351 & 10 & 49.7738 & 15 & 65.1125\\
\end{tabular}
\caption{The first fifteen zeros of $\zeta$ rounded to four decimal places.}\label{table:zeros}
\end{table}

The problem of calculating non-trivial zeros
can be divided into three challenges.
First of all given $s\ne 1$, one wants to be able to compute
$\zeta(s)$ with prescribed precision. Secondly one wants to
locate the critical zeros, the zeros on the half line, up to
a prescribed height $T$. Finally, one wants to show that RH
holds true up to a prescribed height $T$.
 The latter challenge necessitates being able to calculate the total number of zeros in the
critical strip up to height $T$ and comparing this with the
number of critical zeros found.\\
\indent These three challenges are addressed in the following sections.

\section{Euler-Maclaurin formula}
\label{mac}

The first important result on computing non-trivial zeros dates back to
1903, when 
Gram\footnote{J{\o}rgen Pedersen Gram (1850 Nustrup, Denmark -- 1916 Copenhagen, Denmark), Danish mathematician, actuary.}
published a paper with a list of approximate values of 15 non-trivial zeros
and the proof that RH is true up to height $50$ \cite{G}. 
He used the \textbf{Euler-Maclaurin summation formula}, a method 
originating in the 18th century
to evaluate sums, which describes how good an approximation of a sum is obtained
by replacing it by the corresponding integral 
\cite[Chapter 4]{MR2219954}.
A few years before the publication of his paper, he tried a more elaborate method to
compute non-trivial zeros; however, he gave up due to the sheer complexity of the calculations
hoping that someone else would discover a better way\footnote{\textit{Ces difficult\'{e}s m'ayant paru insummontables \`{a} moins de calculs immenses, 
j'abandonnai ces recherches en esp\'{e}rant qu'un autre trouverait quelque m\'{e}thode pouvant servir soit au calcul des coefficients de $\xi(t)$ soit au calcul direct des racines $\alpha$. 
Mais, autant que je sache, aucune m\'{e}thode de ce genre n'a encore \'{e}t\'{e} publi\'{e}e.} \cite{G}}.
Nobody managed to obtain the desired results and after 8 years Gram resumed the work on the numerical estimations of non-trivial zeros.
He tried a "na\" ive approach", which, to his surprise, worked well\footnote{\textit{N\'{e}amoins l'automne dernier je me suis d\'{e}cid\'{e} \`{a} faire cet essai, 
et j'ai \'{e}t\'{e} frapp\'{e} de la facilit\'{e} avec laquelle il a r\'{e}ussi.} \cite{G}}.
Gram applied the Euler-Maclaurin summation formula to $\sum_{j=n}^\infty j^{-s}$ and, evaluating it at $s=1/2+it,$ obtained the (exact), but not absolutely converging, series expansion
\begin{equation}\label{eq:EulerMaclaurin}
 \zeta(s) = \sum_{j=1}^{n-1}j^{-s}+\frac{n^{-s}}{2}+\frac{n^{1-s}}{s-1}+\sum_{k=1}^{\infty}R_{k,n}(s),
\end{equation}
with 
\begin{equation*}
 R_{k,n}(s) = \frac{B_{2k}}{(2k)!}n^{1-s-2k}\prod_{j=0}^{2k-2}(s+j),
\end{equation*}
where $B_k$ denotes the $k$-th Bernoulli number\footnote{Bernoulli numbers can be defined recursively: $B_0=1$ and $B_m = -\sum_{k=0}^{m-1}\binom{m}{k}\frac{B_k}{m-k+1}$, 
or as coefficients of the generating exponential function $\frac{x}{e^{x}-1}=\sum_{n=0}^\infty \frac{B_n x^n}{n!}$.}. 
In order to approximate $\zeta(s)$ one chooses an appropriate $m$, computes the $R_{k,n}$ terms up to this $m$ and estimates the remainder by
\begin{equation*}
\left|\sum_{k=m+1}^{\infty}R_{k,n}(s)\right| < \left|\frac{s+2m+1}{\re s+2m+1}R_{m+1,n}(s)\right|.
\end{equation*}

Gram observed that, in order to calculate any zero with $0<\im s <50$ to 7 decimal places, we can take $n=20$.
In general, the number $n$ of terms needed to be calculated in \eqref{eq:EulerMaclaurin} in order to obtain an estimation of $\zeta(s)$ up to reasonable precision grows linearly with $\im s$.
Later Backlund\footnote{Ralf Josef Backlund (1888 Pietarsaari, Finland -- 1949 Helsinki, Finland), Finnish mathematician, actuary, student of Lindel\"{o}f; 
shortly after obtaining his PhD, he quit academia to work for an insurance company; he was
one of the founders of the Actuarial Society of Finland.} 
provided good bounds for the size of the estimation error for this method in terms of $n$, $m$ and $s$.

This approach is simple and the computational work is proportional to $\im s$. 
However, in practice, it is only efficient to evaluate zeros with small imaginary part.
In Sections \ref{RSformula}, 
\ref{tuurlijk} and \ref{schoon} we will present more efficient algorithms to approximate $\zeta(s)$.

\section{Gram test}

Now that we are able to compute $\zeta(s)$ with arbitrary precision, the next challenge is to locate Riemann's zeros. 
Consider the \textbf{Hardy $Z$-function}\footnote{The first to work
with this function was Riemann, not Hardy.}
\begin{equation}
\label{def:Z}
 Z(t) = e^{i\theta(t)}\zeta\left(\frac{1}{2}+it\right),
\end{equation}
where
\begin{equation}\label{def:theta}
\theta(t)=\textnormal{Im} \log \Gamma\left(\frac{it}{2}-\frac{3}{4}\right) - \frac{t}{2}\log\pi.
\end{equation}
(Recall that $\Gamma$ denotes the Gamma function.)
The function $Z(t)$ is a real-valued function and $Z(t)=0$ if and only if $\zeta\left(\frac{1}{2}+it\right)=0$.
Thus we can detect critical zeros by finding intervals where $Z$ changes sign and then applying Newton's method. 
In particular, we can compute $\theta(t)$ and $\zeta\left(\frac{1}{2}+it\right)$ and then look at the sign change. 
On using well-known asymptotic estimates for the $\Gamma$-function
 we find the asymptotic formula
\begin{equation*}
\theta(t)=\frac{t}{2}\log\left(\frac{t}{2\pi e}\right)-\frac{\pi}{8}+\frac{1}{48t}+O\left(\frac{1}{t^{3}}\right).
\end{equation*}

Alternatively, we could consider $\re \zeta$ and $\im \zeta$ separately. We first observe that $\zeta\left(1/2+it\right) = Z(t)\cos\theta(t)-iZ(t)\sin\theta(t)$.
Then we note that $\im\zeta$ changes sign if either $Z(t)$ or $\sin\theta(t)$ changes sign.
The function $\theta(t)$ is increasing for $t \geq 10$, therefore between consecutive zeros of $\sin\theta(t)$
there is exactly one zero of $\cos\theta(t)$ and so $\cos\theta(t)$ changes sign. 
It follows that, if $\re\zeta(1/2+it)=Z(t)\cos\theta(t)$ is positive at two consecutive zeros $z_1<z_2$ of $\sin\theta(t)$, then $Z$ must change 
sign in the interval $I=[1/2+iz_1,1/2+iz_2]$ 
and hence there is a zero of $Z$ (and hence a 
critical zero) in the interval $I$. 
Thus it is important to locate those $t$ for which $\sin\theta(t)=0$. These points are called \textbf{Gram points}. To be precise, the $n$-th Gram point $g_n$ is defined as the unique solution to
$\theta(g_n)=n\pi$, $g_n \geq 10$. We conclude that, as long as $\zeta\left(1/2+ig_{n}\right)>0$, there is at least one 
critical zero in the interval 
$(1/2+ig_{n-1},1/2+ig_{n})$. Checking if $\zeta\left(1/2+ig_{n}\right)>0$ at a Gram point $g_n$ is called the \textbf{Gram test}.

In 1925
Hutchinson\footnote{John Irwin Hutchinson (1867 Bangor, the US -- 1935), American mathematician, student of Bolza, professor at Cornell University.} \cite{Hu},
using Gram point computations, showed 
that RH holds up to height $T=300$.
He found two empty Gram intervals, but since the
neighbouring Gram interval in each case contained
two zeros, he could overcome this defect.

\section{Total number of zeros up to a given height}
Suppose that we have localized the zeros on the half line for all $0<|\im(s)| < T$. We would like to prove the
RH up to height $T$.
Already Gram showed that 
$$
\log\xi\left(\frac{1}{2}+it\right) = 
\log\xi\left(\frac{1}{2}\right)-
t^2 \sum_{\re \gamma>0}\gamma^{-2}- \frac{t^4}{2} \sum_{\re \gamma>0}\gamma^{-4}-\ldots,
$$
where 
the sum is over the zeros $\rho$ of $\xi$ written in the form
$1/2+i\gamma$. 
From this identity he could evaluate $\sum_{\re \gamma >0}\gamma^{-2n}$ very precisely.
Since $\sum_{\re \gamma >0}\gamma^{-2n}$ is dominated by the first few terms,
he compared the estimated term $\sum \gamma^{-10}$ to the finite sum evaluated for the the 15 zeros found on the critical line.
He then concluded that RH holds up to height $T=50$.
For more zeros the calculations became too involved.

Let $N(T)$ and $N_0(T)$ denote the number of zeros, respectively 
critical zeros, of $\zeta(s)$ with $0<\im s <T$ counted with multiplicities. 
Von Mangoldt (1905) showed that 
\begin{equation}
\label{nntt}
N(T)=\frac{T}{2\pi}\log \frac{T}{2\pi e}+O(\log T).
\end{equation}
\indent Put
\begin{equation}
\label{esthee}
S(T)=\frac{1}{\pi}\im \int_{C_{\epsilon}}\frac{\zeta'(s)}{\zeta(s)}ds,
\end{equation}
where $C_{\epsilon}$ is the path\footnote{It is only a
quarter rectangle and this is related to the four
fold symmetry of the zeta roots.}
consisting of the segment from $1+\varepsilon$ to $1+\varepsilon +iT$ and that from $1+\varepsilon +iT$ to $1/2+iT$. 
Backlund, basing himself on Riemann's ideas, showed in 1912 that
\begin{equation*}
N(T)=\frac{\theta(T)}{\pi}+1+S(T).
\end{equation*}
Moreover, if $\re \zeta$ is not zero on $C_{\epsilon},$  then $N(T)$ is the nearest integer to $\theta(T)/\pi+1$, i.e.
$$|S(T)|\le \frac{1}{2}.$$
Backlund \cite{Ba2, Ba1} himself proved that $N(200)=79$ and, using Gram's approach, showed that all these zeros lie on the half line. 
Hutchinson also used this method to show that there are no more zeros with $0<\im(s)<300$ other than the ones he found.
The problem with the approach of Backlund is to find a $T$ 
such that $\re \zeta$ is non-zero on 
$C$. This is a difficult task. Even worse is the fact that $S(T)$ can get arbitrarily large, so that infinitely often one will pick a $T$ for which $|S(T)|>1/2$. For these two 
reasons, nowadays a 
better method developed by Turing half a century later is used
(see Section \ref{tuurlijk}).

We note that the above approach to verify the RH numerically up
to a certain height will fail if zeta has zeros that are not simple. It is conjectured
that actually the Riemann zeta function has only simple zeros. The reader might
regard this as wishful thinking, but Random Matrix Theory, in particular 
Montgomery's pair correlation conjecture, see Section 
\ref{rrmmtt}, strongly
suggests that all the zeros are simple. It is also a numerical observation that
the critical zeros `repel' each other.\\
\indent We would like to conclude this section by
recalling some major results involving $N(T)$, see 
Karatsuba \cite[Chapter 2]{MR1398651} for an introduction.
Hardy \cite{H} was the first to show that there are infinitely many 
critical zeros. Even more, in 1921 he showed with Littlewood that 
$N_0(T)\ge c_1T$ with $c_1$ a positive constant and $T\ge 15$. In 1942 Selberg
showed that a positive fraction of all non-trivial zeros are critical, that
is that $N_0(T)\ge c_2N(T)$ with $c_2>0$. Selberg's ideas lead
to a tiny value of $c_2$ of about $7 \cdot 10^{-8}$ (as was worked out by 
S. Min in his PhD thesis). In the mid 1970's Levinson\footnote{Norman Levinson (1912 Lynn, the US -- 1975 Boston, the US), American mathematician, professor at MIT.} 
caused a sensation by showing
that one can take $c_2=0.3474$. It was later improved to $0.4088$ by Conrey \cite{MR1004130}, $0.4105$ by Bui, Conrey and Young \cite{MR2825573}, and $0.4128$ by Feng \cite{MR2887604}. 

\section{Riemann-Siegel formula}
\label{RSformula}
In 1932 Siegel \cite{S} published the results on the Riemann zeta function found in the notes of Riemann stored in the archives of the G\"{o}ttingen University Library. 
The following identity was discovered by Riemann. We start by showing the validity
of the integral representation
\begin{equation}\label{eq:ZetaIntegral}
 \zeta(s)=\frac{\Gamma(-s)}{2\pi i} \int_{-\infty}^{+\infty}\frac{(-x)^s}{e^x-1} \frac{dx}{x},
\end{equation}
where the contour integral starts at $+\infty$ and descends along the real axis, circulates $0$ in the positive direction and returns to $+\infty$ up to the real axis. 
We can then evaluate the integral in \eqref{eq:ZetaIntegral} in two ways: firstly by
using the 
trivial identity 
$$\frac{e^{-Nx}}{e^x-1}=\sum_{n=N+1}^\infty e^{-nx},$$ and secondly by changing the contour and applying the Cauchy residue theorem.
We set $\re s=1/2$ and arrive at
\begin{equation*}
Z(t)=2 \sum_{n^2 < t/(2\pi)} \frac{1}{n^{1/2}} \cos(\theta(t)-t\log n) + R(t),
\end{equation*}
where $Z(t)$ and $\theta(t)$ were defined in \eqref{def:Z} 
and \eqref{def:theta} respectively.

The idea of the Riemann-Siegel method\footnote{For a derivation and more
details see, e.g., \cite[Chapter 7]{E}.}
is to bound the remainder term $R(t)$. 
In order to determine the roots of the zeta function, we again use the function $Z$.
Riemann, using the saddle point method, found that 
\begin{equation*}
R(t) = (-1)^{\lfloor \sqrt{t/(2\pi)} \rfloor-1} \left( \frac{t}{2\pi}\right)^{-1/4} \Big(\sum_{j=0}^m C_j t^{-j/2}+R_m(t)\Big),
\end{equation*}
with
$$
 C_0 = F(\delta),~C_1 = -\frac{ F^{(3)}(\delta)}{2^3\cdot 3},~C_2 = \frac{F^{(2)}(\delta)}{2^4} +\frac{F^{(6)}(\delta)}{2^7\cdot 3^2},
$$
etc., and
$$
 F(x) = \frac{\cos\left(x^2+3\pi/8\right)}{\cos\left(\sqrt{2\pi}x\right)},~\delta  =\sqrt{t}+\left(\left\lfloor \sqrt{\frac{t}{2\pi}} \right\rfloor+\frac{1}{2}\right)\sqrt{2\pi}.$$
 The error term $R_m(t)$ satisfies the estimate $R_m(t)=O(t^{-(2m+3)/4})$.
Explicitly estimating the error term $R_m(t)$ is unfortunately quite difficult. 
Titchmarsh\footnote{Edward Charles Titchmarsh (1899 Newbury, England -- 1963 Oxford, England), British mathematician, student of Hardy, professor at the University of Oxford.}, 
a former student of Hardy, 
proved in 1935 \cite{Ti} that if $t > 250\pi,$
then 
\begin{equation*}
 |R_0(t)| \leq \frac{3}{2} \left(\frac{t}{2\pi}\right)^{-3/4}
\end{equation*}
and used this estimate to verify RH up to height 390.

\section{Turing's method}
\label{tuurlijk}
After Siegel's publication in 1932 of the Riemann-Siegel formula, Titchmarsh obtained a
grant for large scale computation of the non-trivial zeros. 
Making use of this formula, tabulating machines and ``computers'' (as the mostly female operators
of such machines were called in those days), Titchmarsh established that the $1\,041$
non-trivial zeros up to height $1\,468$ are all critical.
Turing\footnote{Alan Mathison Turing (1912 London, England -- 1954 Wilmslow, England), British pioneering computer scientist, mathematician, student of Church.}, who from his student days 
onwards was interested in the zeros of the Riemann zeta function, became interested in extending
Titchmarsh's results. 
In his first paper on this topic \cite{Tu1}, he considered the Riemann-Siegel method, but improved the estimate of Titchmarsh's remainder term in such a way that it gives satisfactory results
for $30<\im s<1\,000$. 
However, the paper was rather technical, and, within some years, better estimates were obtained. 
Turing designed and even started to build 
a special purpose analog computer\footnote{The blueprint of the machine is available at \url{http://www.turingarchive.org/browse.php/C/2}.} in order to verify the RH in the range
$0 < \im s< 6\,000$. The outbreak of the Second World War prevented him from completing the construction of the machine.
However, his attempts to build such a machine had set his mind in motion about
building computers and helped him later to build faster his famous Bombes in order
to break the Enigma Code. One can thus argue that many people thank their lives 
to the critical zeros! 

After the Second World War Turing returned to the problem of computing non-trivial zeros.
In his second paper on the zeros of $\zeta$, he attempted to calculate the zeros of the Riemann zeta function
for large $\im s$, that is, for $1\,414<\im s<1\,540,$ in order to find a counterexample. He writes:
\textit{The calculations were done in an optimistic hope that a zero would be found off the critical line, and the calculations were directed more towards finding such zeros than proving that none existed} \cite{Tu2}.
During this attempt to disprove the RH he developed what is now called \textbf{Turing's method}, 
a method to estimate the number of non-trivial zeros up to a given height, which is still used today. 
In order to describe the method, we need to go back to 1914, when Littlewood showed that, 
assuming the RH, the difference
$\pi(x)-{\rm Li}(x)$ changes sign infinitely often(see also 
Section \ref{pivli}), thus
extending a result of Erhard Schmidt who had proved it, assuming the RH was false \cite{L}.
Turing managed to exploit the result of Littlewood in order to average the value of
$$S(t)=N(t)-\Big(\frac{\theta(t)}{\pi}+1\Big)$$
over 
$[0,T]$.
He showed that $S(t)$ tends to zero as $t$ tends to infinity. In particular, he proved that, for $h>0$ and $T>168\pi,$
\begin{equation*}
 \left|\int_{T}^{T+h} S(t) dt\right| \leq 2.3 + 0.128 \log \frac{T+h}{2\pi}.
\end{equation*}
If we calculate $S(t)$ using the observed number of zeros $N'(t)$ instead of $N(t)$ and we missed $k$ zeros, i.e., $N'(t) <N(t)$, then we would get that $S(t) \to -k$, 
which eventually contradicts the above estimate. What the method amounts to is the following. In order to verify RH up to
height $T$, one has to find 
enough critical zeros up to height $T+O(\log T)$.
For more details on Turing's work on this subject, see for example \cite{B}, \cite[265--279]{CL}.

\section{Odlyzko-Sch\"{o}nhage method}
\label{schoon}
In 1988 Odlyzko\footnote{Andrew Michael Odlyzko (1949 Tarn\'{o}w, Poland), Polish-American mathematician, student of Stark, professor at the University of Minnesota.}
and 
Sch\"{o}nhage\footnote{Arnold Sch\"{o}nhage (1934 Lockhausen, the Free State of Lippe), German mathematician, computer scientist, student of Hoheisel, professor at the University of Bonn, 
the University of T\"{u}bingen and the University of Konstanz.} developed a different approach to evaluate the Riemann zeta function. 
In \cite{OS} they observed that the problem of evaluating sums of the form $\sum_{k=1}^Md_kk^{-it}$ 
at an evenly spaced set of $t$'s can be transformed into a problem of evaluating a rational function of the form
$\sum_{k=1}^n a_k (z-b_k)^{-1}$ at all $n$-th roots of unity using the fast Fourier transform. 
The algorithm they designed rapidly evaluates such functions at multiple values. 
In particular, they showed that, for any positive constants $\delta$, $\sigma$, $c_1$, 
there exists an effectively computable constant $c_2$ and an algorithm which, 
for any $T,$ computes $\zeta(\sigma+it)$, $T\leq t\leq T + T^{1/2}$ to $T^{-c_1}$ precision
in less than $c_2T^{\delta}$ operations (addition, subtraction, multiplication, division).
The algorithm has no advantage over the Riemann-Siegel method when a single value is evaluated, 
but it significantly improves the verification time of the RH
for the first $N$ zeros
 as shown
in Table~\ref{table:algorithms}. 

\begin{table}
\footnotesize
 \begin{tabular}{ c | c |  c | c }
  Algorithm & Euler-Maclaurin & Riemann-Siegel  & Odlyzko-Sch\"{o}nhage  \\
  \hline
A & $O(N^{2+\varepsilon})$ & $O(N^{3/2+\varepsilon})$  &  $O(N^{1+\varepsilon})$ \\
  \hline
B & $O(x^{2/3+\varepsilon})$ & $O(x^{3/5+\varepsilon})$ &  $O(x^{1/2+\varepsilon})$ \\
  \hline
C  & $O(x^{\varepsilon})$  & $O(x^{\varepsilon})$ & $O(x^{1/4+\varepsilon})$  \\

\end{tabular}
\caption[caption]{Comparison of different algorithms.\\\hspace{\textwidth}
A: Number of operations needed to verify the RH for the first $N$ zeros;\\\hspace{\textwidth}
B: Number of operations needed to compute $\pi(x)$;\\\hspace{\textwidth}
C: Number of bits of storage needed to compute $\pi(x)$.}\label{table:algorithms}
\end{table}

\section{A brief history of non-trivial zeros calculations}
Riemann seems to have
computed only a few non-trivial zeros. He certainly 
found zeros at approximately 
$1/2+i14.1386$ and at $1/2+i25.31,$ and very likely computed more. He 
derived and tried to
use the wonderful identity
$$\sum_{\im\rho>0}\Big(\frac{1}{\rho}+\frac{1}{1-\rho}\Big)=1+\frac{\gamma}{2}-
\frac{\log \pi}{2}-\log 2=0.02309570896612103381\ldots $$ to prove that the root at 
14.1 is the first root\footnote{The identity can be derived using the
Hadamard factorization for $\xi$. Here $\gamma$ denotes Euler's constant.}. The
latter identity can be used to infer that $\re\rho>10$ for the first non-trivial zero. 
The decimals Riemann gave for the two roots are slightly off, but the 20 decimals above
he computed correctly!\\
\indent In the earlier sections we reported about the computations and innovations due
to Gram, Backlund, Hutchinson, Titchmarsh and Turing. After Turing computers with every
increasing performance played a major role. At the beginning of the 21st century, large scale computations were performed. 
Between 2001 and 2005 the project ZetaGrid, led by Wedeniewski, was established. 
It involved distributed computation of the non-trivial zeros.
It ran on more than $10\,000$ computers in over 70 countries and was
based on software developed by van de Lune, te Riele and Winter
\cite{MR829637} who after years of work eventually computed the first $1.5\cdot 10^9$ Riemann
zeros.
More than $9 \cdot 10^{11}$ first zeros were verified to lie on the critical line.
In 2004 Demichel and Gourdon, \cite{Go}, performed the calculation of the zeros using the Odlyzko-Sch\"{o}nhage method 
to verify that the first $10^{13}$ lie on the critical line. The calculation was not repeated though.
Unfortunately, neither of the results was published in a mathematical journal. \\
\indent The above history is summarized in Table \ref{table:records}.

\begin{table}
\footnotesize
 \begin{tabular}{ c | c |  c }
  Year & Number of computed zeros & Author  \\
  \hline
1859 & $\ge 2$ & Riemann\\  
1903 & 15 & Gram \\
1914 & 79 & Backlund \\
1925 & 138 & Hutchinson \\
1935 & $1\,041$ & Titchmarsh \\
1953 & $1\,104$ & Turing \\
1956 & $25\,000$ & Lehmer \\
1966 & $250\,000$ & Lehman \\
1968 & $3\,500\,000$ & Rosser, Yohe, Schoenfeld \\
1979 & $81\,000\,001$ & Brent \\
1982 & $200\,000\,001$ & Brent, van de Lune, te Riele, Winter \\
1983 & $300\,000\,001$ & van de Lune, te Riele \\
1986 & $1\,500\,000\,001$ & van de Lune, te Riele, Winter \\
2004 & $900\,000\,000\,000$ & Wedeniwski \\
2004 & $10\,000\,000\,000\,000$ & Demichel, Gourdon \\

\end{tabular}
\caption[caption]{Record history of calculating critical  zeros.}\label{table:records}
\end{table}

\section{Applications}
In this section we give some applications of being able to 
determine many 
non-trivial zeros with high precision, e.g., the determination of $\pi(x)$ for 
very large $x$. The three subsections that follow are about prime number inequalities;
the first one about inequalities that hold true, and the other two about famous inequalities conjectured in the 19th century
that turn out to be too good to be true. The final subsection is about the ternary Goldbach
conjecture, which was recently resolved using very extensive zero calculations of
Dirichlet\footnote{Johann Peter Gustav Lejeune Dirichlet, German mathematician (1805--1859), who gave his name to the Dirichlet series; the Dirichlet theorem on primes in arithmetic progressions precedes by almost 70 years the prime number theorem. He was the advisor of Bernhard Riemann.}

$L$-functions (that behave like $\zeta(s)$ in many respects).

\subsection{The exact value of $\pi(x)$}
The record values for
which $\pi(x)$ has been exactly computed are a good indicator for the progress in computational
prime number theory. The obvious
way of computing $\pi(x)$ is, of course, by counting all
primes $p\le x$. For large values of $x$ this is quite
inefficient. In 1871 Meissel devised an ingenious
method of computing $\pi(x)$ without computing all
primes $p\le x$. This method requires only the
knowledge of the primes $p\le \sqrt{x}$, as well
as the values of $\pi(y)$ for some values of
$y\le x^{2/3}$. In 1885 Meissel determined
$\pi(10^9)$ (albeit not quite accurately). The algorithm
was steadily improved and, e.g., in 2007 Oliveira e Silva used
it to compute $\pi(10^{23})$. Assuming RH in 2010 B\"uthe, Franke, Jost
and Kleinjung announced a value of $\pi(10^{24})$.
Very recently by a different method Platt \cite{MR3315519}, based on
an explicit formula of Riemann for the quantity on the left hand side
of (\ref{rielambda}), managed to 
show unconditionally that $\pi(10^{24})=18\,435\,599\,767\,349\,200\,867\,866,$ in
agreement with the value of B\"uthe et al. This computation rests
on the first $69\,778\,732\,700$ 
critical zeros computed with 25 decimal accuracy.
For the values of $\pi(10^k)$ with $8\le k\le 18$ the reader is referred
to Table~\ref{table:pix}.

\subsection{Explicit prime number bounds}
\label{expli}
One of the most often quoted papers in computational number theory is the one by Rosser and Schoenfeld \cite{MR0137689}.
In this paper, among other things, they prove the explicit bounds of $\pi(x)$ 
that we mentioned in the introduction. 
The importance of this paper is that their proof was based
on verifying the RH up to a certain height.
Then using the fact that the first $3\,502\,500$ zeros of $\zeta$ 
lie on the critical line, the authors together with Yohe found 
explicit bounds for the $\vartheta$-function \cite{RYS}.
Let $p_n$ denote the $n$th prime.
Rosser and Schoenfeld obtained the estimates
$$n(\log n+\log \log n-3/2)<p_n<n(\log n+\log \log n-1/2),$$
valid for every $n\ge 21$. {}From this we deduce that
$p_n>n\log n$ for every $n\ge 1$, a result that had been obtained
already in 1939 by Rosser \cite{Ross39}.
Note that Table~\ref{table:pix} suggests that $\pi(x)\ge [x/\log x]$ for all $x$
large enough. Indeed, Rosser and Schoenfeld in 1962 using 
a delicate analysis established the
truth of this inequality for $x\ge 17$.
Meanwhile most explicit estimates of Rosser and Schoenfeld have
been sharpened, e.g., by Dusart who exploited the 
verification of the RH for the first $1.5 \cdot 10^9$ zeros \cite{MR1697455}.

A recent result from B\"uthe shows explicitly the connection between a partial
RH and a sharp prime number estimate. He
proved that if the RH holds true
up to height $T$, then the estimate 
(\ref{schoen1969}) holds true for all $x$ satisfying
$4.92 x/ \log x \le T$.

\subsection{Comparison of $\pi(x)$ with ${\rm Li}(x)$}
\label{pivli}
{}From Riemann's formulae (\ref{R1}) and (\ref{R2}) it follows
that
\begin{equation*}
\pi(x)={\rm Li}(x)-\frac{{\rm Li}(\sqrt{x})}{2}-\sum_{\rho}{\rm Li}(x^{\rho})
+W(x),
\end{equation*}
where $W(x)$ relatively to the three earlier summands is of lower order. 
Riemann's writings suggest (he was rather vague about it) that he 
thought that the inequality 
\begin{equation}
\label{pvl}
\pi(x)<{\rm Li}(x)
\end{equation}
is always satisfied. Gauss and Goldschmidt had established
the validity up to $x=10^5$. Today we know that it even holds up to 
$x=10^{14}$, cf. also Table~\ref{table:pix}. However, in 1914 Littlewood proved the spectacular result
that the difference $\pi(x)-{\rm Li}(x)$ changes sign infinitely often.
In the mid-1930s Ingham showed that this result follows from 
knowledge of some initial non-trivial zeros\footnote{Subsequently, Turing whose initial interest 
was in improving the results of Skewes, turned his interest to 
computing non-trivial zeros.}. His proof was both simpler and
more explicit than Littlewood's, but also more computational. 
Let $x_0$ be the smallest integer for which $\pi(x_0)>{\rm Li}(x_0)$.
Skewes showed
that
$$x_0<10^{10^{10^{34}}} {\rm{(1933,~on~RH),}} \quad x_0<10^{10^{10^{964}}} {\rm(1955,~unconditionally)}.$$
For a long time these bounds of Skewes were considered to be the largest ``naturally'' occurring 
numbers in mathematics.\\
\indent Using tables of non-trivial zeros accurate to 28 digits for the first $15\,000$ zeros and to 14 digits for the next $35\,000$ zeros, te 
Riele showed that (\ref{pvl}) is false for at least $10\,180$ successive integers in $[6.627\cdot 10^{370},6.687\cdot 10^{370}]$. He made use of an earlier
result of Lehman, which allows one to put bounds on 
$\pi(x)-{\rm Li}(x)$, assuming that we have found all 
non-trivial zeros to height $T$, and that RH is checked up to height $T$
\cite{MR0202686}. Lehman himself had obtained $x_0<10^{1166}$ using this result.

We now know that once (\ref{pvl}) holds the wait until the inequality
is reversed grows again, on average, as a function of the starting $x$. 
Thus we should perhaps not be surprised that the average of $\chi(t)$
with $\chi(t)=1$ if $\pi(t)<{\rm Li}(t)$ and $\chi(t)=0$ otherwise, does not exist.
However, under various plausible conjectures the density
\begin{equation*}
\delta=\lim_{x\rightarrow \infty}\frac{1}{\log x}\int_1^x \frac{\chi(t)}{t}dt
\end{equation*}
does exist and satisfies $\delta=0.99999973\ldots$ (see \cite{MR1329368}). This conjecturally
quantifies the dominance of Li$(x)$ over $\pi(x)$. The strong bias towards Li$(x)$ is an example
of an interesting phenomenon that is called Chebyshev's bias (see Section \ref{zz}).

\subsection{The Mertens conjecture}
\label{merti}
The \textbf{Mertens function} $M(x)$ denotes the difference between the number of squarefree
integers $n\le x$ having an even number of prime factors and those
having an odd number of prime factors. 
More formally, we have 
$$M(x)=\sum_{n\le x}\mu(n),$$ where
$\mu(n)=0$ if a square exceeding one divides $n$ and $\mu(n)=(-1)^m$ with $m$ being 
the number of different prime factors of $n$ otherwise. The function 
$\mu$ is called the 
\textbf{M\"obius\footnote{August Ferdinand M\"{o}bius (1790 Schulpforta, Saxony -- 1868 Leipzig, Germany), German mathematician and astronomer, student of Gauss and Pfaff, 
professor at the University of Leipzig.} 
function} and outside number theory arises very frequently in combinatorial counting
(\textbf{inclusion-exclusion}).
It is not difficult to show
that 
$$\frac{1}{\zeta(s)}=\sum_{n=1}^{\infty}\frac{\mu(n)}{n^s}=s\int_1^{\infty}\frac{M(t)}{t^{1+s}}dt,$$
from which one can infer that the RH is equivalent with
$|M(t)|=O(t^{1/2+\epsilon})$. 
The PNT is equivalent to $M(x)=o(x)$.
Most number theorists envision the Mertens function as something like
a random walk, with the $\pm 1$ contributions from $\mu$ to $M(x)$
being close to a random coin flip. 
It is known by the probability theory that for the summatory function
$w(x)=\sum_{n\le x}w_n$, with $w_n=\pm 1$, randomly and independently, we
have
$$\lim \sup_{x\rightarrow \infty} \frac{w(x)}{\sqrt{(x/2)\log \log x}}=1.$$
This suggests that if $\mu$ is "sufficiently"
random, then the ratio $M(x)/\sqrt{x}$ is expected to be unbounded.

On the basis of numerical work 
Stieltjes\footnote{Thomas Joannes Stieltjes (1856 Zwolle, the Netherlands -- 1894 Toulouse, France), Dutch mathematician, professor at Toulouse University.} 
in 1885\footnote{In a 
letter to Hermite published only in 1905.} and
independently Mertens\footnote{Franz Carl Joseph Mertens (1840 Schroda, Prussia -- 1927 in Vienna, Austria), Polish-German mathematician, student of Kummer and Kronecker, 
professor at Jagiellonian University in Cracow and the University of Vienna.}
 in 1897 \cite{M}, conjectured that
$|M(x)|<\sqrt{x}$, a conjecture that is now known as the \textbf{Mertens conjecture}.
Daublebsky von Sterneck\footnote{
Robert Daublebsky von Sterneck (1871 Vienna, Austria--1928 Graz,Austria), a student of F. Mertens. He was born into a celebrated family belonging to the nobility. Later in life he was poor and only used Sterneck as a last name.},
earlier had made the stronger claim that
$|M(x)|<\sqrt{x}/2$ for $x>200$. That conjecture was disproved in 1963 by Neubauer. We now
know that $7\,725\,038\,629$ is the minimal integer $>200$ for which the Daublebsky von Sterneck
bound does not hold.
In 1983 Odlyzko and te Riele \cite{MR0783538} caused a sensation by disproving the Mertens conjecture.
They showed that
$$\lim \sup_{x\rightarrow \infty} \frac{M(x)}{\sqrt{x}} >1.06,\quad \lim \inf_{x\rightarrow \infty}\frac{M(x)}{\sqrt{x}}<-1.009,$$
without actually giving a specific integer $x_0$ with $|M(x_0)|\ge \sqrt{x_0}$.
Later
Kotnik and te Riele
\cite{MR2282922} using very extensive computer calculations showed 
that there is an $x_0<e^{1.59\cdot 10^{40}}$ for which the Mertens conjecture
fails.\\
\indent The first step in disproving 
the Mertens conjecture is to find the analogue of (\ref{psipsi2}) for 
$M(x)/\sqrt{x}$. Assuming RH and that all critical zeros are simple this can be done. 
It involves a sum having 
terms of the form
$e^{i\gamma y}/(\rho \zeta'(\rho))$, where $\rho=1/2+i\gamma$,  
$x=e^y$. 
Next, one finds an upper bound for the error
made on cutting this formula at a given
height $T$ up to which one has
verified RH numerically and computed the zeros with sufficient
numerical precision.
In order to obtain the desired result, one needs to verify RH and 
the simplicity of the zeros up to a height large enough. Secondly, we want
each of the terms to be close to its maximum, which happens when $\gamma y$ 
is close to an even integer. The latter leads to a problem in simultaneous
Diophantine approximation. Now the key factor that allowed Odlyzko and 
te Riele to progress beyond earlier failed attempts to disprove the Mertens
conjecture, was using the at the time recently developed Lenstra-Lenstra-Lov\'asz 
algorithm, or \textbf{LLL-algorithm} for short. With the use of this new
algorithm less extensive computing was needed to reach rather stronger results 
than before. Indeed only the first $2\,000$ 
non-trivial zeros 
calculated with circa 100 significant decimal places, 
were used in the disproof.

\subsection{The ternary Goldbach conjecture}
\label{gold}
The \textbf{Goldbach conjecture} (formulated in 1742 in a letter to Euler) states 
that every even number exceeding 2 can be represented as a sum of two primes. 
This is a very famous conjecture that remains unsolved. A weaker variant is
the \textbf{ternary Goldbach conjecture}, also known as odd or weak Goldbach conjecture.
It says that every odd number greater than 7 can be expressed as the sum of three odd
primes.

Hardy and Littlewood showed in 1923 that on 
GRH\footnote{Grand Riemann Hypothesis, 
see Section \ref{zz}.}, the odd Goldbach conjecture is true for all sufficiently large odd numbers. In 1937, 
Vinogradov \cite{Vin37} 
established this result unconditionally. 
Vinogradov used the \textbf{circle method}, which involves both minor
and major arc estimates.
His student Borozdkin in his 
PhD thesis (unpublished),
made the sufficiently large explicit, yielding a huge number $e^{e^{e^{41.96}}}$ and further published a result with the 
smaller bound $e^{e^{16.038}}$, see \cite{Bor56}.
In 1997, Deshouillers, Effinger, te Riele and Zinoviev published a result showing that GRH implies Goldbach's weak conjecture \cite{MR1469323}. 
This required checking all integers $\le 10^{20}$ as for all
larger integers they could establish the result by
theoretical means.
Without GRH it was known in 2002 that the $10^{20}$ has
to be replaced by $10^{1347}$.
In 2012 and 2013, Helfgott released a series of preprints improving 
the major and minor arc estimates sufficiently to unconditionally prove the weak Goldbach conjecture (see \cite{Helb}, \cite{Helc}, \cite{Held} and \cite{MR3201598}).
Helfgott made use of a result of Platt \cite{MR3522979} who had rigorously
verified, using interval arithmetic, the RH for all Dirichlet
$L$-functions for modulus $q\le 400\,000$ up to height around
$10^8/q$.
The binary Goldbach conjecture is numerically verified up to $4 \cdot 10^{18}$ by Oliveira e Silva, Herzog and Pardi \cite{MR3194140}.\\
\indent The Goldbach conjecture is also studied in the field of \textbf{additive number theory},
where one considers a special set $A$ (primes, squares, etc.) and then wonders what the
sumset $A+A=\{a+b: a,b \in A\}$ looks like (similarly, with $A+A+A,$ etc.). For a nice introduction
to additive number theory see the book by Tao and Vu \cite{MR2573797}.

\section{Major recent developments}
In this section we discuss some major more recent developments related 
to the Riemann zeta function.
\subsection{Complex analytic number theory}
Riemann's paper and, in its wake, the proof of the PNT were a major achievement of 19th century
mathematics and gave rise to complex analytic number theory.\\
\indent A function $f$ from the natural numbers to the complex numbers is called
an \textbf{arithmetic function}. An important subclass are the \textbf{multiplicative
functions} that satisfy $f(1)=1$ and $f(mn)=f(m)f(n)$ for arbitrary coprime natural
numbers $m$ and $n$ (the M\"obius function is an example of a multiplicative function).
The behaviour of arithmetic functions is usually very erratic. For that
reason it makes sense to investigate related quantities that show more
regular behaviour, for example the summatory function $\sum_{n\le x}f(n)$. The zeta
function reflex leads one to consider $F(s)=\sum_{n=1}^{\infty}f(n)n^{-s}$, which is
called a \textbf{Dirichlet series}.
Using the Perron integral and assuming that $F(s)$ converges absolutely for some $\re s>\sigma$ one then finds that 
$$\sum_{n\le x}f(n)= \int_{c-i\infty}^{c+i\infty} F(s)\frac{x^s}{s}ds,$$ 
with $c>\sigma$ arbitrary and tries to estimate the integral. For this it is particularly important to be able
to shift the line of integration as far as possible to the left. Here it becomes relevant to have an analytic continuation of $F(s)$. The more
one knows about $F(s)$ the more information one can deduce about
$f$. However, often $F(s)$ is not so well-behaved and one tries to
consider a closely related function $f^*$ instead, with the property that its
Dirichlet series $F^*(s)$ behaves in a nicer way\footnote{E.g., 
$\sum p^{-s}$ is not nicely behaved, but $\sum \Lambda(n)n^{-s}$ is.}.
The results one obtains
in this way about $f^*$ one tries to relate back to $f$. \\
\indent The attentive reader sees of course immediately that this whole approach is
patterned on Riemann's 1859 paper.

\subsection{The Zeta Zoo}
\label{zz}
The study of the Riemann zeta-function
has proved so extraordinarily successful
in deepening our understanding of the primes,
that it has become standard in number theory
to try to associate zeta type functions to
arithmetic structures (some kind of Pavlov reflex). Indeed, these days there
is an enormous zoo of zeta functions. The alpha
animal in the Zeta Zoo was and still is zeta. It
codifies the behaviour of the integers and their
atomic constituents: the primes. An important
species of zeta functions are the \textbf{Dirichlet $L$--functions}.
These were introduced by Dirichlet in order to
understand the behaviour of primes in arithmetic
progression. Like zeta, they satisfy a product
expansion, a functional
equation and it is conjectured that their zeros are
also on the half line. The latter hypothesis is
called the \textbf{Extended Riemann Hypothesis} (ERH).
Again a lot of computational number theory was
done to verify RH for individual Dirichlet $L$--series \footnote{A profound database can be found at \url{http://www.lmfdb.org/}.}
up to a certain height and this was used to derive
explicit bounds for 
$\pi(x;d,a),$ with $a$ and $d$ being coprime integers. 
These computations played an important role in
the recent proof of the ternary Goldbach conjecture
by Helfgott (see Section \ref{gold}).
De la Vall\'ee Poussin (1897) proved that every
of the arithmetic progressions $a({\rm mod~}d)$ with $1\le a<d$ and 
$a$ and $d$ coprime (of which there are $\varphi(d)$)
gets asymptotically its fair share of the primes, i.e., that asymptotically
$$\pi(x;d,a)\sim \frac{\pi(x)}{\varphi(d)}.$$
In 1837 Dirichlet in a ground breaking paper 
(where he introduced characters in number theory) had proved a weaker version of this result.
Although the primes are asymptotically 
equidistributed, they have some positive bias towards
progressions modulo $d$ where $a$ is a non-square modulo $d$.
This was noted in 1853 by Chebyshev and is now known as
\textbf{Chebyshev's bias}.
E.g., Bays and Hudson found in 1979 that $608\,981\,813\,029$ is the
smallest prime for which $\pi(x;3,2)>\pi(x;3,1)$. For a nice 
introduction to this phenomenon see Granville and Martin
\cite{MR2158415}. As in the $\pi(x)$ versus Li$(x)$ problem,  under
various assumptions there is a computable logarithmic measure 
for how often $\pi(x;d,a_1)>\pi(x;d,a_2)$.\\
\indent Another important class of zeta functions are the
\textbf{Dedekind zeta functions}. This is the analogue of the
Riemann zeta function for a number field and can be treated 
by similar methods. E.g., Landau in 1903 showed the 
\textbf{prime ideal theorem}, stating that in a given number field
the number of prime ideals of norm $\le x$ grows asymptotically
as $x/\log x$. The hypothesis that every Dedekind zeta function
has its non-trivial zeros on the critical line is called
the \textbf{Grand Riemann Hypothesis} (GRH).\\
\indent The reader might wonder about a more stringent definition of a
zeta function. Here is what Selberg, with a life time of experience with 
zeta functions, thought about this.\\
\indent A zeta function $F(s)$ is a function of a complex variable 
$s$ that satisfies the following properties.
\begin{enumerate}
\item Dirichlet series: for $\re s > 1$, 
one can write $F(s) = \sum_{n=1}^{\infty}a_nn^{-s}$.
\item Ramanujan hypothesis: the growth of the coefficients $a_n$ has to be modest, in essence like that of the divisor function $\sum_{d\mid n}1.$
\item Analytic continuation: $F(s)$ extends to a meromorphic function.
\item Functional equation: there is "a connection" between $F(s)$ and $F(1-s)$.
\item Euler product\footnote{Patterned
after Euler's product formula \eqref{prodie}.} one should be able to write $F(s)=\prod_p F_p(s).$
\end{enumerate}

\subsection{Correlation of pairs of non-trivial zeros}
\begin{figure}
\centering
\includegraphics[width=\linewidth]{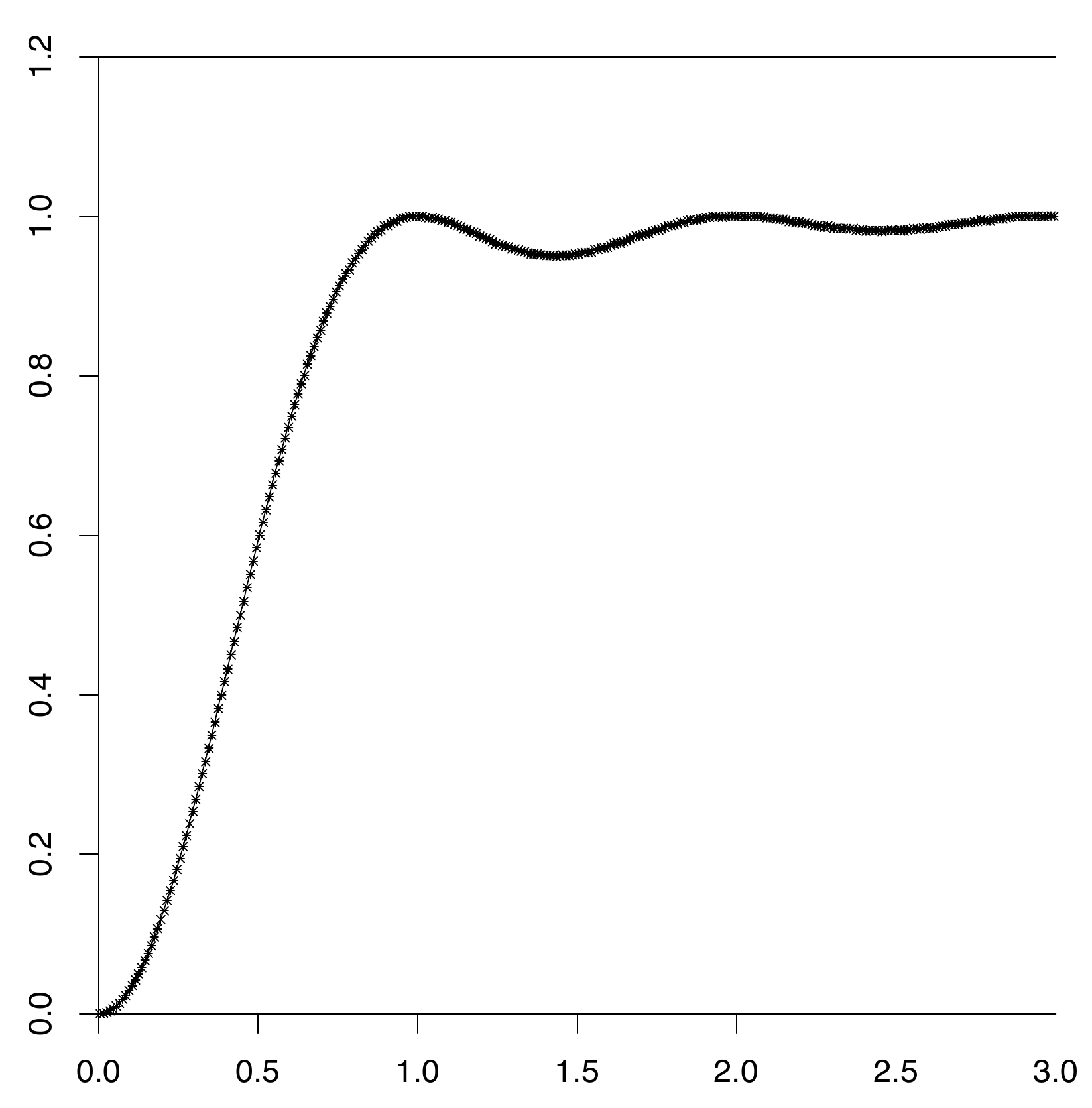}
\caption{Odlyzko's pair correlation plot for $2\cdot 10^8$ non-trivial
zeros near the $10^{23}$th zero. In the displayed interval, the data
agrees with the pair correlation conjecture to within about $0.002$.}
\label{fig:rmt}
\end{figure}

So far we exclusively focused on finding and counting non-trivial zeros. A refinement is
to ask about the distribution of the zeros, e.g., how are the differences between
(consecutive) zeros distributed? Here the first theoretical work is due to Montgomery\footnote{Hugh Lowell Montgomery (1944 Muncie, the US), American mathematician, student of Davenport, professor at Michigan University.} \cite{MR0337821}.
He assumed the RH and wrote the zeros as $1/2+i\gamma_i$ with $0<\gamma_1\le \gamma_2\le \ldots$
To account for the increase of the density of the zeros as one goes up the critical strip
and rescales them as 
$${\overline \gamma_i}=\gamma_i\frac{\log(\gamma_i/(2\pi e))}{2\pi}.$$
By the asymptotic (\ref{nntt}) for
$N(T)$ it then follows that the mean spacing between two 
rescaled consecutive zeros is
1. The \textbf{pair correlation conjecture} of
Montgomery states that, with $0\le \alpha<\beta$, we have, as $M$ tends
to infinity,
$$\frac{1}{M}|\{1\le i<j\le M:{\overline \gamma_j}-{\overline \gamma_i}\in [\alpha,\beta)\}|
\sim \int_{\alpha}^{\beta}\Big(1-\Big(\frac{\sin(\pi t)}{\pi t}\Big)^2\Big)dt.$$
The integrand here is small when $t$ is close to zero, suggesting that the non-trivial zeros
repel each other. Montgomery proved a smoothened version of his conjecture and used it to
show that on the RH more than two thirds of the non-trivial zeros are simple.\\
\indent Odlyzko numerically tested both this conjecture and the RH, beginning in the late
1980s. In the 1990s this led to monumental computations where billions of zeros were computed
high up the critical strip.  For a graphical 
 demonstration of the computations, see Figure~\ref{fig:rmt}.
\subsection{Random matrix theory}\label{rrmmtt}
In 1972 Montgomery
discussed his pair correlation conjecture with the physicist Dyson\footnote{Freeman John Dyson (1923 Crowthorne, England), American physicist, professor at Princeton University.}.
Dyson immediately saw that the statistical distribution found by Montgomery appeared to be the same as the pair correlation distribution for the eigenvalues of a random Hermitian matrix that he had discovered a decade earlier.\\
\indent \textbf{Random matrix theory} (RMT) was proposed by the physicist 
Eugene Wigner\footnote{Eugene Paul Wigner (1902 Budapest, Hungary -- 1995 Princeton, the US), Hungarian-American physicist and mathematician, Laureate of the Nobel Prize in Physics (1963).} 
in 1951 to describe nuclear physics. The quantum mechanics of a heavy nucleus is complex and 
not well understood. Wigner made the bold conjecture that the statistics of 
the energy levels can be captured by random matrices.

RMT turns out to be a powerful tool in making conjectures involving the 
Riemann zeta and related functions. These conjectures lie typically way beyond
the reach of current theoretical tools. Since this is the case, doing numerical
checks on the conjectures is very important. These checks are often very
computationally intensive.

As a very important example let us consider the problem of determining the even moments
of the Riemann zeta function on the half line. We define
$$I_k(T)=\frac{1}{T}\int_0^{T}|\zeta(1/2+it)|^{2k}dt.$$
Progress on determining the moments was slow.
Hardy and Littlewood \cite{MR1555148}
determined in 1916 the asymptotic behaviour of $I_1(T)$.
Ingham \cite{MR1575391} improved this in 1926 in two ways, calculating the full asymptotic expansion when $k=1$ and determining the leading term for $k=2$ (the lower-order terms for 
$k=2$ were determined by Heath-Brown \cite{MR532980} in 1979).
It was conjectured for general $k$ that
$$I_k(T)\sim \frac{f(k)a(k)}{k^2!}T(\log T)^{k^2},$$ with
$a(k)$ a certain infinite product over all primes and $f(k)$ 
an integer.
Hard grinding analytic number theory
led to the conjecture that 
$f(3)=42$ \cite{MR1639551} and $f(4)=24024$ \cite{MR1828303}.
Progress beyond this seemed very challenging. Using RMT and modelling
the zeros up to height $T$ with $N$ by $N$ matrices with
$N$ around $\log(T/2\pi)$, Keating and Snaith came with a conjectural
integer value for $f(k)$ for all $k$. Their conjecture is consistent
with all earlier results and the conjectural values of $f(3)$ and $f(4)$ 
derived using hardcore analytic number theory. Extensive numerical work
corroborates the Keating and Snaith conjecture. 

The RMT method offers a dictionary that allows one to translate number
theoretical problems into random matrix problems that usually can 
be solved. However, there is no proof whatsoever that the dictionary
always works. As long as this is the case, computational number theoretical
work will play a very important role.
The work
of Odlyzko \cite{O} on the zeros near $10^{22}$ is of enormous 
importance here, as 
the Riemann zeta function starts showing
``its true face'' only for very large values of $T$. \\
\indent The RMT method also works in the setting of \textbf{function fields}. These share
many similarities with the number field setting, but often are easier to deal with.
E.g., for them the Riemann Hypothesis is proved! 
In this setting Katz and Sarnak actually managed to prove various important results suggested by 
the RMT dictionary (published in their book \cite{MR1659828} 
and surveyed in \cite{MR1640151}).

\indent The quest for
an explanation of the RMT connection is ongoing and has led to active
research at the intersection of number theory, mathematical physics, probability
and statistics.

\section{Words of warning}
We hope that we have whetted the appetite of the reader to do 
computations in analytic number theory him or herself. A word of warning is, however,
not amiss. Asymptotic estimates in analytic number theory often involve
repeated logarithms. These are very difficult to detect by numerical computation.
Thus a function that grows, e.g., like $\log \log \log x,$ looks on a computer like a bounded function. For this reason it is highly dangerous 
to make conjectures based on numerics 
alone (e.g., the Mertens conjecture), without some theoretical and heuristic 
considerations supporting the truth of the conjecture. 

\section{Further reading}
Nice popular introductions to prime number theory are
Sabbagh \cite{Sabbagh} and du Sautoy \cite{duSautoy}. Halfway between a popular
and more mathematical treatment is an interesting collection 
of prime number records and results modelled after the Guiness
book of records \cite{MR1377060}.
For more 
mathematical introductions see, e.g.,
\cite{C,MR1790423,E,MR1074573,P}. 
In the book of Edwards\footnote{Warning: Edwards 
writes $\Gamma(s+1)$ instead of $\Gamma(s)$.} there is a detailed explanation of the method used
by Gram, Backlund and Hutchinson to compute the smallest 300 non-trivial zeros.
More advanced books are Ivi\'{c} \cite{MR792089} and Titchmarsh 
\cite{MR882550}.
The contribution of Hejhal and Odlyzko \cite[265--279]{CL} on the work
of Turing 
on the Riemann zeta function
is very informative.
Snaith \cite{MR2684776} wrote a nice overview of the connections between $L$-functions and random matrix theory, focusing on the example of the Riemann zeta function. 
A very readable conference proceedings on this matter is \cite{MR2145172}.
Rubinstein's beautiful article 
\cite[pp. 633--679]{MR3525906} also discusses
RMT, but with main focus on the influence of Riemann.\\ 
\indent This survey owes a lot to the book of Narkiewicz \cite{N} on the
development of prime number theory that provides a nice mix of mathematical ideas
and historical material. Also the book of Crandall and Pomerance 
\cite{MR2156291} on computational
prime number theory turned out to be quite helpful.\\
\indent Finally, in tune with the dictum of Edwards "that one should read
the masters and beware of secondary sources"\footnote{Such as the present survey written by non-masters.}, we
would like to point out \cite{MR2463715}, which has many articles by the Riemann zeta masters.

\section{Acknowledgement}
The authors are enormously grateful to Alexandru Ciolan
for his \TeX nical assistance and his excellent proofreading of the many earlier versions. 
Also they thank Alex Weisse for his \TeX nical assistance.
Special thanks are due to Jan B\"uthe, Andreas Decker, 
Tomoko Kitagawa (alias Kate Kattegat),
W\l adys\l aw Narkiewicz, 
Olivier Ramar\'e, Michael Rubinstein, 
Kannan Soundararajan, 
Caroline Turnage-Butterbaugh and Tim Trudgian.

\bibliographystyle{abbrv}
\bibliography{mybib}

\end{document}